\documentclass[12pt]{article}
\usepackage{amsmath,amscd,amsthm,amssymb}

\usepackage{color}

\def\Rn{{\mathbb{R}^n}}

\def\a {\alpha}
\def\i{\infty}

\def\L1loc{L_1^{\rm loc}(\Rn)}

\theoremstyle{definition}

\theoremstyle{remark}

\numberwithin{equation}{section}

\def\a{\alpha}

\def\i{\infty}

\begin{document}
	\theoremstyle{definition}
	\newtheorem{definition}{Definition}[section]
	\newtheorem{example}[definition]{Example}
	\newtheorem{remark}[definition]{Remark}
	\newtheorem{observation}[definition]{Observation}
	\theoremstyle{plain}
	\newtheorem{theorem}[definition]{Theorem}
	\newtheorem{lemma}[definition]{Lemma}
	\newtheorem{proposition}[definition]{Proposition}
	\newtheorem{corollary}[definition]{Corollary}
	\numberwithin{equation}{section}

\begin{center}
{\Large \bf{Maximal and fractional maximal operators in the Lorentz-Morrey spaces and their applications to the Bochner-Riesz and Schr\"{o}dinger-type operators}}
\end{center}

\

\centerline{\large Abdulhamit Kucukaslan$^{a,b,}$\footnote{Corresponding author
		\\
$~~~~~$E-mail address: kucukaslan@pau.edu.tr}}

\

\centerline{$^{a}$\it Institute of Mathematics of Czech Academy of Sciences,}
\centerline{ \it 115 67, Prague, Czech Republic}

\centerline{$^{b}$\it Faculty of Applied Sciences, Pamukkale University, 20680, Denizli, Turkey}

\

\begin{abstract}
The aim of this paper is to obtain boundedness conditions for the maximal function $Mf$ and to prove the necessary and sufficient conditions  for the fractional maximal oparator $M_{\alpha}$ in the Lorentz-Morrey spaces $\mathcal{L}_{p,q;\lambda}(\Rn)$ which are a new class of functions. We get our main results by using the obtained sharp rearrangement estimates. The obtained results are applied to the boundedness of particular operators such as the Bochner-Riesz operator $B_r^\delta$ and the Schr\"{o}dinger-type operators 	$V^{\gamma} (-\Delta+V)^{-\beta}$ and $V^{\gamma} \nabla (-\Delta+V)^{-\beta}$ in the  Lorentz-Morrey spaces $\mathcal{L}_{p,q;\lambda}(\Rn)$, where the nonnegative potential $V$ belongs to the reverse H\"{o}lder class $B_{\infty}(\Rn)$.

\

\noindent{\bf AMS Mathematics Subject Classification:} 42B20, 42B35, 47G10

\noindent{\bf Keywords:} Maximal operator, Fractional maximal operator, Lorentz-Morrey spaces,  Bochner-Riesz operator, Schr\"{o}dinger-type operators.

\end{abstract}

\

\section{Introduction}

Let $B(x,r)$  the open ball centered at $x$ of radius $r$ for $x\in \Rn$ and $|B(x,r)|$ is the Lebesgue measure of $B(x,r)$.  The fractional maximal operator is defined  at $f\in L_{1}^{loc}(\Rn)$ by
$$M_{\alpha}f(x):=\sup_{r>0}|B(x,r)|^{\frac{\alpha}{n}-1}\int_{B(x,r)}|f(y)|dy, 0\leq \alpha<n $$
where the supremum is taken over all the balls centered at $x$ of radius $r$.
Note that in the case $\alpha=0$
we get the classical Hardy-Littlewood maximal operator $M:=M_{0}$. It is well known that for the  maximal operator $M$ the rearrangement inequality
$$
cf^{\ast\ast}(t) \leq (Mf)^{\ast}(t)\leq Cf^{\ast\ast}(t), ~~t \in (0,\infty)
$$
holds, (see \cite{BenSh}, Chapter 3, Theorem 3.8) where $f^{*}$ is the non-increasing rearrangement of $f$ such that $f^{*}(t):=inf \left\{{\lambda > 0: d_{f}(\lambda) \le 0}\right\}$,
$d_{f}(\lambda)$ denotes the distribution function of $f$ given by $d_{f}(\lambda):=| \left\lbrace x \in (0,\infty):\left| f(x)\right|> \lambda \right\rbrace  |$
for all  $t>0$, and $f^{**}(t):=\frac{1}{t}\int_0^t f^{\ast}(s)ds.$

The Lorentz-Morrey spaces $\mathcal{L}_{p,q;\lambda}(\Rn)$ are a new class of functions and introduced by Mingione in \cite{Mg} as follows. 
\begin{definition}
	 Let $ 1\leq p < \infty $, $ 0 < q < \infty $, $ 0\leq \lambda \leq n,$ and   $f \in {\mathcal L}_{p,q;\lambda}(\Rn)$. Then the Lorentz-Morrey space $\mathcal{L}_{p,q;\lambda}(\Rn)$ is the set of all measurable functions $f$ on $\Rn$  iff
	\begin{align*} 
	\|f\|_{{\mathcal L}_{p,q;\lambda}(\Rn)}&:=\sup_{x \in \Rn, \, r>0} r^{-\frac{\lambda}{p}} \|f\chi_{{B(x,r)}}\|_{L_{p,q}(\Rn)}<\infty.
	\end{align*}
\end{definition}
Mingione  \cite{Mg}, studied the boundedness of the restricted fractional maximal operator $M_{\alpha, B_{0}}$ in the restricted Lorentz-Morrey spaces $\mathcal{L}_{p,q;\lambda}(B)$, where $B_{0}$ is a given ball and $B$ is any other ball contained in $B_{0}$ and containing $x$. The author derived a general non-linear version, extending \textit{a priori} estimates and \textit{regularity} results for possibly degenerate non-linear elliptic problems to the various spaces of Lorentz and Lorentz-Morrey type considered in \cite{Adams1982, AGKS, Mg} and \cite{Ragusa2012}.  In \cite{Ragusa2012}, Ragusa studied some embeddings between these spaces. Note that the spaces ${\mathcal L}_{p,q;\lambda}(\Rn)$ and $L_{p,q;\lambda \frac{q}{p}}(\Rn)$ defined by Mingione and Ragusa respectively, coincide, thus
${\mathcal L}_{p,q;\lambda}(\Rn)= L_{p,q;\lambda \frac{q}{p}}(\Rn).$ The local variant of Lorentz-Morrey spaces $\mathcal{L}_{p,q;\lambda}(\Rn)$ replacing by $B(0,r)$ instead of $B(x,r)$, so called the local Morrey-Lorentz spaces $\mathcal{L}_{p,q;\lambda}^{loc}(\Rn)$ are introduced and the basic properties of these spaces are given in \cite{AGS}.  Recently, in \cite{AGKS, GAKS} and \cite{GKAS}, the authors studied the boundedness of some classical operators of harmonic analysis in these spaces.

In this paper,  first, we give some basic properties of Lorentz-Morrey spaces $\mathcal{L}_{p,q;\lambda}(\Rn)$. Furthermore, we get the sharp rearrangement inequalities which we use while proving our results. Next, in section 3,  we obtain the boundedness conditions for the maximal function $Mf$ in the  Lorentz-Morrey spaces $\mathcal{L}_{p,q;\lambda}(\Rn)$ and  we get the necessary and sufficient conditions for boundedness of the fractional maximal operator $M_{\alpha}$ in the spaces $\mathcal{L}_{p,q;\lambda}(\Rn)$. Finally, in section 4, we apply these results to  the Bochner-Riesz operator $B_r^\delta$ and the Schr\"{o}dinger-type operators 	$V^{\gamma} (-\Delta+V)^{-\beta}$ and $V^{\gamma} \nabla (-\Delta+V)^{-\beta}$ in the  Lorentz-Morrey spaces $\mathcal{L}_{p,q;\lambda}(\Rn)$, respectively, where the nonnegative potential $V$ belongs to the reverse H\"{o}lder class $B_{\infty}(\Rn)$.

Throughout the paper, we denote by $c$ and $C$ for positive constants, independent of appropriate parameters and not necessary the same at each occurrence. If $p \in \left[ 1,\infty\right]$, the conjugate number $p'$ is defined by $\frac{1}{p}+\frac{1}{p'}=1$. Finally, for non-negative expressions $A_{1},A_{2}$ we use the symbol $A_{1} \approx A_{2}$ to express that $cA_{1} \leq A_{2}\leq CA_{1}$ for some positive constants $c$ and $C$ independent of the variables in the expressions $A_{1}$ and $A_{2}$.
\section{Preliminaries}
	The Lorentz space $ L_{p,q}(\Rn)$  is the collection of all measurable functions of $f$ on $\Rn$ such the quantity
	\begin{align*}\label{kuc01}
	\left\| f\right\|_{L_{p,q}(\Rn)}:= \Big\{
	\begin{array}{cc}
	\left( \int_{0}^{\i} \left(t^{\frac{1}{p}}f^{\ast}(t) \right)^{q}\frac{dt}{t} \right) ^{\frac{1}{q}}, & 0<q<\i,  0<q<\i \\
	\sup_{t>0} t^{\frac{1}{p}}f^{\ast}(t), & 0< p \leq \i, q=\i,
	\end{array}%
	\end{align*} 
	is finite. If $1<p \leq\i,1\leq q \leq \i$, then  
	
	\begin{equation*} \label{kuc02}
	\|f\| _{L{{p,q}(\Rn)}}\leq \|f\|_{L{{p,q}(\Rn)}}^{* }
	\leq \frac{p}{p-1}\|f\|_{L{{p,q}(\Rn)}}.
	\end{equation*}
	 For more detail useful references about Lorentz spaces considered in \cite{BenSh}.

We denote by $ L_{p,\lambda}(\Rn)$ Morrey space given in \cite{Morrey}; $0\leq \lambda \leq n, 1\leq p \leq \i, f \in L_{p, \lambda} $, if $f \in L_{p}^{loc}(\Rn)$ and 
\begin{equation*}\label{kuc05}
\| f\| _{L_{p,\lambda}(\Rn)}:=\sup \limits_{x \in \Rn, r>0} r^{-\frac{\lambda}{p}} \, \|f\| _{L_{p}(B(x,r))}<\i.
\end{equation*}
Morrey spaces appeared to be useful in the study  of local behavior properties of the solutions of second order elliptic PDEs.  For more information about Morrey-type spaces see \cite{ BurGulHus1, BG1, GulIsmKucSJFS2015, Kuc1} and \cite{Kuc2}.

The Lorentz-Morrey spaces $\mathcal{L}_{p,q;\lambda}(\Rn)$ are a very natural generalization of the Lorentz spaces $L_{p,q}(\Rn)$ and Morrey spaces $L_{p,\lambda}(\Rn)$.
\begin{remark}
As a consequence by Lemma \ref{kuc26} $(ii)$, if $q=p$
 then $\mathcal{L}_{p,p,\lambda}(\Rn)\equiv L_{p,\lambda}(\Rn)$, if $\lambda=0$
 then $\mathcal{L}_{p,q,0}(\Rn) \equiv L_{p,q}(\Rn)$, and $\lambda=n, p=q$, then $\mathcal{L}_{p,p,n}(\Rn)\equiv L_{\i}(\Rn)$. 
If $\lambda < 0$ or $\lambda > n$, then $\mathcal{L}_{p,q; \lambda}(\Rn)\equiv \theta$, where $\theta$ the set of all functions equivalent to $0$ on $\Rn$.
\end{remark}
\begin{lemma}\cite{BenSh}, \cite{EE}, \cite{St1}\label{kuc26}
	
	(i) Let $0<p<\i$, then $\int_{\Rn} \left| f(x)\right|^p dx=\int_0^\i \left( f^{\ast}(t)\right) ^pdt$ holds.
	
	(ii) For any $t>0$, $\sup_{\left| E\right| =t}\int_{E} \left| f(x)\right| dx=\int_0^t f^{\ast}(s)ds.$
	
	(iii)  For any $t>0$, $\left( f+g\right)^{\ast}(t)\leq f^{\ast}\left(\frac{t}{2} \right) + g^{\ast}\left(\frac{t}{2} \right) $ holds.
\end{lemma}
\begin{lemma}\label{kuc25}
	Let $0 \leq \alpha <n$. Then there exist a positive constant $C$, depending on $\alpha$ and $n$ such that
	\begin{equation}
	\label{kuc08}
	\sup\limits_{t>0}t^{1-\frac{\alpha }{n}}\left(  M_{\alpha }f\chi_{B(x,t)}\right) ^{*}( t) \leq C\int\limits_{\Rn}|f(x)|dx
	\end{equation}
	and
	\begin{equation}
	\label{kuc09}
	\sup\limits_{t>0} \left(  M_{\alpha}f\chi_{B(x,t)}\right) ^{*}(t) \leq C\sup\limits_{t>0}t^{\frac{\alpha }{n}
	}f^{*}(t).
	\end{equation}  	
\end{lemma}

\begin{proof}
	The estimate \eqref{kuc08} follows from  (Theorem 1.1, in \cite{Cia}).
	For the estimate \eqref{kuc09}, for every $B(x,r) \subset \Rn$, we get 
	\begin{align*}
	\sup_{r>0}|B(x,r)|^{\frac{\alpha}{n}-1}\int_{B(x,r)}|f(y)|dy&\leq |B(x,r)|^{\frac{\alpha}{n}-1}\int_{0}^{|B(x,r)|}t^{\frac{\alpha}{n}}f^*(t)t^{-\frac{\alpha}{n}}dt
	\\
	&\leq\frac{n}{n-\alpha}\sup\limits_{t>0} t^{\frac{\alpha}{n}}f^*(t).
	\end{align*}
Hence the proof is completed.
\end{proof}

\begin{lemma}\label{kuc16}
	Let $0 \leq \alpha < n$. Then there exist a positive constant $C$, depending only on $n$ and $\alpha$, such that
	\begin{equation}  \label{kuc10}
	\left(M_{\alpha}f\chi_{B(x,t)}\right)^{* }(t)\leq C\sup\limits_{t<\tau <\infty }
	\tau ^{\frac{\alpha }{n }}\left(f\chi_{B(x,r)}\right)^{**}(\tau ),~~t>0
	\end{equation}
	holds for all $f\in L_{1}^{loc}(\Rn)$.
	Inequality $\eqref{kuc10}$ is sharp in the sense that
	for all $\varphi \in \mathcal{M}^{+}( 0,\i ;\downarrow ) $
	there exists a function $f$ on $\Rn$ such that $f^{*}=\varphi $ a.e. on $(0,\i)$ and
	\begin{equation}\label{kuc11}
		\left(M_{\alpha}f\chi_{B(x,t)}\right)^{* }(t)\geq c\sup\limits_{t<\tau <\infty }
	\tau ^{\frac{\alpha }{n }}\left(f\chi_{B(x,r)}\right)^{**}(\tau ),~~t>0,
	\end{equation}
	where $\mathcal{M}^{+}(0,\i ;\downarrow)$ is the set of all non-negative and non-increasing measurable functions on $(0,\i) $ and $c$ is a positive constant which depends only on $n$ and $\alpha$.
\end{lemma}

\begin{proof}
	To prove the inequality \eqref{kuc10}, we may suppose that
	\begin{equation*}\label{kuc12}
	\sup\limits_{	t<\tau <\i }\tau ^{\frac{\alpha }{n}}\left(f\chi_{B(x,r)}\right)^{**}(\tau )<\i,
	\end{equation*}
	otherwise there is nothing to prove. Then by Lemma \ref{kuc26} $(i)$
	\begin{equation*}
	\int_{\Rn} \left| f\chi_{B(x,r)}(x)\right| dx=\int_0^t \left(f\chi_{B(x,r)}\right)^{\ast}(s) ds
	\end{equation*}
	holds for all $E \subset \Rn$ with $\left| E\right| \leq t$. In particular, if we put 
	\begin{equation*}
	E=\left\lbrace x: \left| f(x)\right| > \left(f\chi_{B(x,r)}\right)^{\ast}(t)\right\rbrace 
	\end{equation*}
	then $\left| E\right| \leq t$ and so $f \in L_1(E)$. Then the function 
	\begin{equation*}
	g_{t}(x)=\max \left\lbrace \left| f(x)\right|- \left(f\chi_{B(x,r)}\right)^{\ast}(t), 0 \right\rbrace sgn f(x),
	\end{equation*}
	belongs to $L_1(\Rn)$. Also the function
	\begin{equation*}
	h_{t}(x)=\min \left\lbrace \left| f(x)\right|,\left(f\chi_{B(x,r)}\right)^{\ast}(t) \right\rbrace sgn f(x),
	\end{equation*}
	holds 
	\begin{equation*}
	(h_{t})^*(\tau)=\min \left\lbrace \left(f\chi_{B(x,r)}\right)^{\ast}(\tau),\left(f\chi_{B(x,r)}\right)^{\ast}(t) \right\rbrace, \tau \in (0, \i).
	\end{equation*}
	Thus
	\begin{align}\label{kuc13}
	\sup\limits_{\tau>0 }\tau ^{\frac{\alpha }{n}}(h_{t})^*(\tau)&=\max \left\lbrace \sup\limits_{
		0<\tau <t } t ^{\frac{\alpha }{n}} \left(f\chi_{B(x,r)}\right)^{\ast}(t), \sup\limits_{
		t\leq\tau <\i } \tau ^{\frac{\alpha }{n}} \left(f\chi_{B(x,r)}\right)^{\ast}(\tau) \right\rbrace \notag
	\\
	&=\sup\limits_{
		t\leq\tau <\i } \tau ^{\frac{\alpha }{n}} \left(f\chi_{B(x,r)}\right)^{\ast}(\tau)\leq \sup\limits_{
		t\leq\tau <\i } \tau ^{\frac{\alpha }{n}} \left(f\chi_{B(x,r)}\right)^{**}(t)
	\end{align}
	which together with  the inequality \eqref{kuc11} implies that $h_{t} \in WL_{\frac{n}{\alpha}}$. Furthermore, since $f=h_{t}+g_{t}$ and 
	\begin{equation}\label{kuc14}
	(g_{t})^*(\tau)=\chi_{[0,t)}(\tau)\left(  \left(f\chi_{B(x,r)}\right)^{\ast}(\tau)-\left(f\chi_{B(x,r)}\right)^{\ast}(t) \right) , \tau \in (0, \i).
	\end{equation}
	By using Lemma \ref{kuc26} $(iv)$, Lemma \ref{kuc25}, the inequalities \eqref{kuc13} and \eqref{kuc14}, we get
	\begin{align*}
	\left(M_{\alpha}f\right)^{*}(t)& \leq \left(M_{\alpha}g_{t}\right)^{*}\left(\frac{t}{2} \right)+\left(M_{\alpha}h_{t}\right)^{*}\left(\frac{t}{2} \right)
	\\
	&\lesssim \left(\left(\frac{t}{2} \right)^{\frac{\alpha}{n}-1} \int_{\Rn}g_{t}(y)dy+ \sup\limits_{\tau>0 }\tau ^{\frac{\alpha }{n}}(h_{t})^*(\tau)\right) 
	\\
	&\lesssim t^{\frac{\alpha}{n}-1} \int_{0}^t \left(  \left(f\chi_{B(x,r)}\right)^{\ast}(\tau)-\left(f\chi_{B(x,r)}\right)^{\ast}(t) \right) d \tau
	\\
	&+ \sup\limits_{0<\tau<\i }\tau ^{\frac{\alpha }{n}}(f\chi_{B(x,r)})^{**}(\tau) 
	\\
	&\lesssim \sup\limits_{0<\tau<\i }\tau ^{\frac{\alpha }{n}}(f\chi_{B(x,r)})^{**}(\tau)
	\end{align*}
	and the inequality \eqref{kuc10} follows.
	Furthermore, the inequality \eqref{kuc11} exist for all $t \in (0,\i)$. Let $\varphi \in \mathcal{M}^{+}(0,\i ;\downarrow)$, where $\mathcal{M}^{+}(0,\i ;\downarrow)$ is the set of all non-negative and non-increasing measurable functions on $(0,\i) $. Putting $f(x)=\varphi(\omega_{n}\left|x \right|^{n} )$, where $\omega_{n}$ is the volume of the unit ball in $\Rn$, $\omega_{n}=\left|B(0,r) \right| $; and $y \in B(x,r)$, we have $(f\chi_{B(x,r)})^*=\varphi(0,\i)$. Moreover, denote by $B(x,\left| y\right| )$ the ball with centered $x$ and having radius $\left| y\right| $. Then, for $\left| y\right| >\left| x\right| $,
	\begin{align*}
	\left(M_{\alpha}f\right)^{*}(t)& = \sup_{r>0}|B(x,r)|^{\frac{\alpha}{n}-1}\int_{B(x,r)}|f(y)|dy
	\\
	&\geq |B(x,\left| y\right|)|^{\frac{\alpha}{n}-1}\int_{B(x,\left| y\right|)}|f(y)|dy 
	\\
	&= c\left( \omega_{n} (\left| y\right|^{n})^{\frac{\alpha}{n}-1}\int_{0}^{\omega_{n}\left| y\right|^{n}} (f\chi_{B(x,r)})^{*}(\tau)d\tau \right)
	\\
&=cH(\omega_{n} \left| y\right|^{n}) ,
	\end{align*}
	where $H$ is the Hardy operator given in \cite{Sam} defined as $H(t)=t^{\frac{\alpha}{n}-1}\int_{0}^{t}\varphi(\tau)d\tau, t \in (0,\i)$. Consequently, 
	\begin{align*}
	\left(M_{\alpha}f\right)(x) \geq c\sup_{\tau >\omega_{n} \left| x\right|^{n}}H(\tau),
	\end{align*}
	thus  the inequality \eqref{kuc11} follows on taking rearrrangements. Hence the proof is completed.
\end{proof}
\section{Main Results}

In this section, we characterize the boundedness conditions of maximal operators $M$ and prove the necessary and sufficient conditions  for the fractional maximal oparator $M_{\alpha}$ in the Lorentz-Morrey spaces $\mathcal{L}_{p,q;\lambda}(\Rn)$ by using the obtained sharp rearrangement estimates.

\begin{theorem}\label{kuc17}
Let $1< p < \i, 1 \leq  q < \i$ and  $0\leq \lambda \leq n$ and for all  $f \in \mathcal{L}_{p,q; \lambda}(\Rn)$, then maximal operator $M$ is bounded  in the Lorentz-Morrey spaces $\mathcal{L}_{p,q;\lambda}(\Rn)$.
\end{theorem}

\begin{proof}
Let $1< p < \i, 1 \leq  q < \i$. Then by using definition of the spaces $\mathcal{L}_{p,q; \lambda}(\Rn)$, Lemma \ref{kuc26} $(ii)$ and Lemma \ref{kuc25}, we get 
	\begin{align*}
	\left\| Mf\right\| _{\mathcal{L}_{p,q; \lambda}(\Rn)}&=\sup_{r>0}r^{-\frac{\lambda}{p}} \left\|t^{\frac{1}{p}-\frac{1}{q}} \left(  M f\right)  ^{* }(t)\right\| _{L_{q}(0,\i)}
	\\
	&=\sup_{r>0}r^{-\frac{\lambda}{p}}\left(\int_{0}^{\i}\left(t^{\frac{1}{p}}\left(  M f\right)  ^{* }(t) \right)^{q}\dfrac{dt}{t}  \right)^{\frac{1}{q}} 
	\\
	&\approx \sup_{r>0}r^{-\frac{\lambda}{p}}\left(\int_{0}^{\i}\left(t^{\frac{1}{p}}\left(  f\chi_{B(x,r)}\right)  ^{** }(t) \right)^{q}\dfrac{dt}{t}  \right)^{\frac{1}{q}}
		\\
	&\leq \sup_{r>0}r^{-\frac{\lambda}{p}} \left\|f \right\|^{*}_{L_{p,q}(B(x,r))} 
	\\
	&\leq \frac{p}{p-1} \left\| f\right\| _{\mathcal{L}_{p,q; \lambda}(\Rn)}.
	\end{align*} 
	Hence the maximal operator $M$ is bounded on the Lorentz-Morrey spaces $\mathcal{L}_{p,q; \lambda}(\Rn)$.
\end{proof}

\begin{theorem}\label{kuc18}
	Let $0 \leq \alpha < n$. Then the following statements are equivalent:
	
	$(i)$ If $1< p \leq  q < \i, 1 \leq u \leq s \leq \i, 1<p<\frac{n-\lambda}{\alpha}, 0<\lambda<n$, then  the fractional maximal operator $M_{\alpha}$ is bounded from  Lorentz-Morrey space $\mathcal{L}_{p,u,\lambda}(\Rn)$ to another one $\mathcal{L}_{q,s,\lambda}(\Rn)$ such that
	\begin{align*}
	\left\| M_{\alpha }f\right\| _{\mathcal{L}_{q,s; \lambda}(\Rn)}\lesssim \left\| f\right\| _{\mathcal{L}_{p,u; \lambda}(\Rn)}.
	\end{align*}
	
	$(ii)$ For all $\varphi \in \mathcal{M}^{+}( 0,\i ;\downarrow ) $ 
	there exists a positive constant $C$ such that
	\begin{align}\label{kuc23}
	\sup_{r>0}r^{-\frac{\lambda}{q}} &\left[\int_{0}^{\i} \left( \sup_{t<\tau<\i}\tau^{\frac{\alpha}{n}-1}\int_{0}^{\tau} \varphi(\sigma)d\sigma\right) ^{s}t^{\frac{s}{q}-1}dt \right]^{\frac{1}{s}}\notag
	\\
	&\leq C \sup_{r>0}r^{-\frac{\lambda}{p}}\left[\int_{0}^{\i}\varphi^{p}(t)t^{\frac{u}{p}-1}dt \right]^{\frac{1}{u}}  .
	\end{align} 
	
	$(iii)$ $ \frac{1}{p}-\frac{1}{q}=\frac{\alpha}{n-\lambda} $.
\end{theorem}
	
\begin{proof}
$(i) \Leftrightarrow (ii)$.

$(i)$ Assume that the fractional maximal operator $M_{\alpha}$ is bounded from $\mathcal{L}_{p,u,\lambda}(\Rn)$ to $\mathcal{L}_{q,s,\lambda}(\Rn)$. Then 
$ \left\| M_{\alpha }f\right\| _{\mathcal{L}_{q,s,\lambda}(\Rn)}\lesssim \left\| f\right\| _{\mathcal{L}_{p,u,\lambda}(\Rn)}$
holds.

For every $ \varphi=\left( f\chi_{B(x,r)}\right)  ^{* }(t) \in \mathcal{M}^{+}( 0,\i ;\downarrow ) $, $\left( f\chi_{B(x,r)}\right)  ^{* }=\varphi$ a.e. on $(0,\i)$, and from Lemma \ref{kuc16}
\begin{align*}
\sup_{r>0}r^{-\frac{\lambda}{q}} &\left[\int_{0}^{\i} \left( \sup_{t<\tau<\i}\tau^{\frac{\alpha}{n}-1}\int_{0}^{\tau} \left( f\chi_{B(x,r)}\right)  ^{* }(\sigma)d\sigma\right) ^{s}t^{\frac{s}{q}-1}dt \right]^{\frac{1}{s}}
\\
&= \sup_{r>0}r^{-\frac{\lambda}{q}} \left[\int_{0}^{\i} \left( \sup_{t<\tau<\i}\tau^{\frac{\alpha}{n}} \left( f\chi_{B(x,r)}\right)  ^{** }(\tau)\right) ^{s}t^{\frac{s}{q}-1}dt \right]^{\frac{1}{s}}  
\\
&\lesssim \sup_{r>0}r^{-\frac{\lambda}{q}} \left[\int_{0}^{\i} \left(  \left(M_{\alpha } f\right)  ^{* }(t)\right) ^{s}t^{\frac{s}{q}-1}dt \right]^{\frac{1}{s}} 
\\
&\lesssim \sup_{r>0}r^{-\frac{\lambda}{p}} \left[\int_{0}^{\i} \left(  \left( f\chi_{B(x,r)}\right)  ^{* }(t)\right) ^{u}t^{\frac{u}{p}-1}dt \right]^{\frac{1}{u}}
\end{align*} 
holds.

$(ii)$ Conversely,  for every $ \varphi=\left( f\chi_{B(x,r)}\right)  ^{* }(t) \in \mathcal{M}^{+}( 0,\i ;\downarrow ) $, $\left( f\chi_{B(x,r)}\right)  ^{* }=\varphi$ a.e. on $(0,\i)$, and from Lemma \ref{kuc16}
\begin{align*}
\left\| M_{\alpha }f\right\| _{\mathcal{L}_{q,s,\lambda}(\Rn)}&= \sup_{r>0}r^{-\frac{\lambda}{q}}\left[\int_{0}^{\i} \left(  \left(M_{\alpha } f \right)^{* }(t)\right) ^{s}t^{\frac{s}{q}-1}dt \right]^{\frac{1}{s}} 
\\
&\lesssim  \sup_{r>0}r^{-\frac{\lambda}{q}}  \left[\int_{0}^{\i} \left( \sup_{t<\tau<\i}\tau^{\frac{\alpha}{n}} \left( f\chi_{B(x,r)}\right) ^{** }(\tau)\right) ^{s}t^{\frac{s}{q}-1}dt \right]^{\frac{1}{s}}
\\
&= C   \sup_{r>0}r^{-\frac{\lambda}{q}} \left[\int_{0}^{\i} \left( \sup_{t<\tau<\i}\tau^{\frac{\alpha}{n}-1} \int_{0}^{\tau} \left( f\chi_{B(x,r)}\right)  ^{* }(\sigma)d\sigma \right) ^{s}t^{\frac{s}{q}-1}dt \right]^{\frac{1}{s}}  
\\
& \lesssim  \sup_{r>0}r^{-\frac{\lambda}{p}} \left[\int_{0}^{\i} \left(  \left( f\chi_{B(x,r)}\right)  ^{* }(t)\right) ^{p}t^{\frac{u}{p}-1}dt \right]^{\frac{1}{u}}
\\
&=\left\| f\right\| _{\mathcal{L}_{p,u,\lambda}(\Rn)}
\end{align*}
holds.

$(ii) \Leftrightarrow (iii)$
The equivalence of $ (ii) $ and $ (iii) $ follows from the same proof method in \cite{Opic}.
Hence the proof is completed.
\end{proof}

\section{Some Applications}
\subsection{The estimate of Bochner-Riesz operator in the spaces $\mathcal{L}_{p,q; \lambda}(\Rn)$}

Let $\delta>(n-1)/2$,
$B_r^\delta(f)^{\hat{}}(\xi)=(1-r^2|\xi|^2)_+^\delta\hat{f}(\xi)$
and $B_r^\delta(x)=r^{-n}B^\delta(x/r)$ for $r>0$. The maximal
Bochner-Riesz operator is defined by (see \cite{LiuLu} and \cite{LiuChen})
$$
B_{\delta,*}(f)(x)=\sup_{r>0}|B_{r}^\delta(f)(x)|.
$$
It is clear that (see \cite{GarRub})
\begin{align}\label{kuc444}
	B_{\delta,*}(f)(x) \lesssim Mf(x).
\end{align}
Since the maximal operator $M$ is bounded on the  Lorentz-Morrey spaces $\mathcal{L}_{p,q; \lambda}(\Rn)$, then from Theorem \ref{kuc17} we get the following statement.

\begin{theorem} \label{ak13} Let $1< p < \i, 1 \leq  q < \i$ and  $0\leq \lambda \leq n$, and there exist a positive constant $C$ independent of $f$ and for all  $f \in \mathcal{L}_{p,q; \lambda}(\Rn)$. Then the Bochner-Riesz operator $B_r^\delta$ is bounded on the  Lorentz-Morrey spaces $\mathcal{L}_{p,q; \lambda}(\Rn)$.
\end{theorem}
\begin{proof}
	The idea of proofs of Theorem \ref{ak13} is based on the inequality \eqref{kuc444} in which the maximal Bochner-Riesz operator $B_{\delta,*}$  dominated by the operator $M$. Hence, the proof is step by step the same as in the proof of Theorem \ref{kuc17}. 
\end{proof}
For the case $\lambda=0,$ from Theorem \ref{ak13} we get the following statement.
\begin{corollary} \label{kuc02} Let $1< p < \i, 1 \leq  q < \i$ and  $0\leq \lambda \leq n$. Then the Bochner-Riesz operator $B_r^\delta$ is bounded on the  Lorentz spaces $L_{p,q}(\Rn)$.
\end{corollary}

\subsection{The estimates of Schr\"{o}dinger-type operators 	$V^{\gamma} (-\Delta+V)^{-\beta}$ and $V^{\gamma} \nabla (-\Delta+V)^{-\beta}$ in the spaces $\mathcal{L}_{p,q; \lambda}(\Rn)$}

When $V$ is a non-negative polynomial, Zhong (\cite{Zhong}) proved that the operators $V^{k}(-\Delta + V)^{-k}$
and $V^{k-1/2}\nabla (-\Delta + V)^{-k}$, $k\in \mathrm{N}$, are bounded on $L_{p}(\Rn)$, $1<p\leq\infty$.
Shen \cite{Shen} studied the Schr\"{o}dinger operator $-\Delta + V$, assuming the nonnegative potential $V$
belongs to the reverse H\"{o}lder class $B_q(\Rn)$ for $q \ge n/2$ and he proved the $L_p(\Rn)$ boundedness of the operators
$(-\Delta + V)^{i \gamma}$, $\nabla^{2} (-\Delta+V)^{-1}$, $\nabla (-\Delta+V)^{-\frac{1}{2}}$ and $\nabla (-\Delta+V)^{-1}$.

We give the boundedness of the Schr\"{o}dinger-type operators
$$
\mathcal{T}_1=V^{\gamma} (-\Delta+V)^{-\beta}, ~~0 \le \gamma \le \beta \le 1,
$$
and
$$
\mathcal{T}_2=V^{\gamma} \nabla (-\Delta+V)^{-\beta},
~~ 0 \le \gamma \le \frac{1}{2} \le \beta \le 1,
~ \beta-\gamma \ge \frac{1}{2}
$$
 from the Lorentz-Morrey spaces $\mathcal{L}_{p,u;\lambda}(\Rn)$ to another one $ \mathcal{L}_{q,s;\lambda}(\Rn)$.
Note that the operators $V (-\Delta+V)^{-1}$ and $V^{\frac{1}{2}} \nabla (-\Delta+V)^{-1}$
in \cite{Li} are the special case of $\mathcal{T}_1$ and $\mathcal{T}_2$, respectively.

It is worth pointing out that we need to establish pointwise estimates for
$\mathcal{T}_1$, $\mathcal{T}_2$ by using the
estimates of fundamental solution for the Schr\"{o}dinger
operator on $\Rn$ in \cite{Li}. Then we prove the boundedness of the Schr\"{o}dinger-type operators 	$V^{\gamma} (-\Delta+V)^{-\beta}$ and $V^{\gamma} \nabla (-\Delta+V)^{-\beta}$ in the Lorentz-Morrey spaces $\mathcal{L}_{p,q; \lambda}(\Rn)$ by using boundedness of the fractional maximal operators $ M_{\alpha } $ in these spaces.

The following two pointwise estimates for $\mathcal{T}_1$ and $\mathcal{T}_2$ are proved
in \cite{Sugano}  with the potential $V \in B_{\infty}$.

\noindent{\bf Theorem A.}  {\it \cite{Sugano} Suppose that $V \in B_{\infty}$ and $0 \le \gamma \le \beta \le 1$.
	Then for any $f \in C_{0}^{\infty}(\Rn)$
	$$
	\left|\mathcal{T}_1 f(x)\right| \lesssim  M_{\a}f(x),
	$$
	where $\a=2(\beta-\gamma)$.
}

\noindent{\bf Theorem B.} {\it \cite{Sugano}  Suppose that $V \in B_{\infty}$,
	$0 \le \gamma \le \frac{1}{2} \le \beta \le 1$ and $\beta-\gamma \ge \frac{1}{2}$.
	Then for any $f \in C_{0}^{\infty}(\Rn)$
	$$
	\left|\mathcal{T}_2 f(x)\right| \lesssim M_{\a}f(x),
	$$
	where $\a=2(\beta-\gamma)-1$.
}

From Theorem \ref{kuc18} and by using Theorems A and B we get the following two statements, respectively.
\begin{theorem}  \label{kuc21}
	Let $V \in B_{\infty}$, $0 \le \gamma \le \beta \le 1$. Then the following statements are equivalent:
	 
$(i)$ If $1< p \leq  q < \i, 1 \leq u \leq s \leq \i, 1<p<\frac{n-\lambda}{2(\beta-\gamma)}, 0<\lambda<n$.	Then the Schr\"{o}dinger-type operator $\mathcal{T}_1$ is bounded from $\mathcal{L}_{p,u;\lambda}(\Rn)$ to $\mathcal{L}_{q,s;\lambda}(\Rn)$, such that there is a positive constant $C$ the inequality
\begin{align*}
\|\mathcal{T}_1 f\|_{\mathcal{L}_{q,s;\lambda}(\Rn)} \le C \|f\|_{\mathcal{L}_{p,u;\lambda}(\Rn)}
\end{align*}
holds for all $f \in C_{0}^{\infty}(\Rn) \cap \mathcal{L}_{p,q; \lambda}(\Rn)$. 

$(ii)$ The inequality \eqref{kuc23} holds for all $\varphi \in \mathcal{M}^{+}( 0,\i ;\downarrow ) $.

$(iii)$ $ \frac{1}{p}-\frac{1}{q}=\frac{2(\beta-\gamma)}{n-\lambda} $.	
\end{theorem}

\begin{theorem}  \label{kuc22}
Let $V \in B_{\infty}$, $0 \le \gamma \le \frac{1}{2} \leq \beta \le 1, \beta - \gamma \geq \frac{1}{2}$.  Then the following statements are equivalent: 

$(i)$ If $1< p \leq  q < \i, 1 \leq u \leq s \leq \i, 1<p<\frac{n-\lambda}{2(\beta-\gamma)-1}, 0<\lambda<n$.	Then the Schr\"{o}dinger-type operator $\mathcal{T}_2$ is bounded from $\mathcal{L}_{p,u;\lambda}(\Rn)$ to $\mathcal{L}_{q,s;\lambda}(\Rn)$, such that there is a positive constant $C$ the inequality
	\begin{align*}
	\|\mathcal{T}_2 f\|_{\mathcal{L}_{q,s;\lambda}(\Rn)} \le C \|f\|_{\mathcal{L}_{p,u;\lambda}(\Rn)}
	\end{align*}
		holds for all $f \in C_{0}^{\infty}(\Rn) \cap \mathcal{L}_{p,q; \lambda}(\Rn)$.
		
$(ii)$ The inequality \eqref{kuc23} holds for all $\varphi \in \mathcal{M}^{+}( 0,\i ;\downarrow ) $.
	
$(iii)$	 	$ \frac{1}{p}-\frac{1}{q}=\frac{2(\beta-\gamma)-1}{n-\lambda} $. 

\end{theorem}

\begin{proof}
	The idea of proofs of Theorem \ref{kuc21} and Theorem \ref{kuc22} are based on the Theorem A and Theorem B in which the Schr\"{o}dinger-type operators $\mathcal{T}_1$ and $\mathcal{T}_2$ dominated by the operator $M_{\alpha}$, respectively. Hence, the proofs are step by step the same as in the proof of Theorem \ref{kuc18}. 
\end{proof}

\end{document}